\documentclass[letterpaper, 10pt, twocolumn, conference]{ieeeconf}
\IEEEoverridecommandlockouts
\usepackage{amssymb, amsmath, mathrsfs, subfigure, url}
\usepackage[dvipsnames, usenames]{color}
\usepackage{graphicx}

\let\oldsqrt\sqrt
\def\sqrt{\mathpalette\DHLhksqrt}
\def\DHLhksqrt#1#2{%
\setbox0=\hbox{$#1\oldsqrt{#2\,}$}\dimen0=\ht0
\advance\dimen0-0.2\ht0
\setbox2=\hbox{\vrule height\ht0 depth -\dimen0}%
{\box0\lower0.4pt\box2}}

\DeclareMathOperator{\erf}{erf}
\DeclareMathOperator{\erfc}{erfc}

\DeclareMathOperator{\diag}{diag}
\DeclareMathOperator{\sat}{sat}
\DeclareMathOperator{\sgn}{sgn}

\def\NN{\mathbb{N}}
\def\Nz{\mathbb{N}_0}
\def\R{\mathbb{R}}
\def\posR{\mathbb{R}_{\ge 0}}
\def\EE{\mathbb{E}}

\def\x0{x_0}
\def\xz{x_0}

\newcommand{\iss}{\textsc{iss}}
\newcommand{\ClassKL}{\mathcal{KL}}
\newcommand{\ClassKinfty}{\mathcal{K}_\infty}
\newcommand{\fa}{\forall\,}

\renewcommand{\subset}{\subseteq}
\renewcommand{\le}{\leqslant}
\renewcommand{\ge}{\geqslant}

\newcommand{\tr}[1]{\mathbf{tr}\!\left(#1\right)}
\newcommand{\nn}{\nonumber}

\newcommand{\transp}{^\mathsf{T}}
\newcommand{\mrm}{\mathrm}
\newcommand{\smat}[1]{\left[\begin{matrix} #1 \end{matrix}\right]}
\newcommand{\abs}[1]{\left\lvert{#1}\right\rvert}

\newcommand{\RemEnd}{\hspace{\stretch{1}}{$\vartriangleleft$}}
\newcommand{\ExEnd}{\hspace{\stretch{1}}{$\triangle$}}
\newcommand{\DefEnd}{\hspace{\stretch{1}}{$\Diamond$}}

\newcommand{\lambdamax}[1]{\lambda_{\text{max}}(#1)}

\newcommand{\norm}[1]{\left\lVert{#1}\right\rVert}
\newcommand{\epower}[1]{\mathrm e^{#1}}

\newcommand{\indic}[1]{\mathbf 1_{#1}}

\newtheorem{theorem}{Theorem}

\newtheorem{lemma}[theorem]{Lemma}
\newtheorem{proposition}[theorem]{Proposition}
\newtheorem{example}[theorem]{Example}

\newtheorem{defn}[theorem]{Definition}

\newtheorem{remark}[theorem]{Remark}

\title{On Stochastic Model Predictive Control with Bounded Control Inputs\thanks{This research was partially supported by the Swiss National Science Foundation under grant 200021-122072}}
\author{Peter Hokayem, Debasish Chatterjee, John Lygeros
\thanks{The authors are with the Automatic Control Laboratory, Electrical Engineering, ETH Zurich, Switzerland \texttt{hokayem,chatterjee,lygeros@control.ee.ethz.ch}}
}

\begin{document}
\maketitle

\begin{abstract}
This paper is concerned with the problem of Model Predictive Control and Rolling Horizon Control of discrete-time systems subject to possibly unbounded random noise inputs, while satisfying hard bounds on the control inputs. We use a nonlinear feedback policy with respect to noise measurements and show that the resulting mathematical program has a tractable convex solution in both cases. Moreover, under the assumption that the zero-input and zero-noise system is asymptotically stable, we show that the variance of the state, under the resulting Model Predictive Control and Rolling Horizon Control policies, is bounded. Finally, we provide some numerical examples on how certain matrices in the underlying mathematical program can be calculated off-line.
\end{abstract}

\section{Introduction}
Model Predictive Control (MPC) for deterministic systems has received a considerable amount of attention over the last few decades, and significant advancements have been realized in terms of theoretical analysis as well as industrial applications. The motivation for such research thrust comes primarily from tractability of calculating optimal control laws for constrained systems. In contrast, the counterpart of this development for stochastic systems is still in its infancy.

The deterministic setting is dominated by worst-case analysis relying on robust control methods. The central idea is to synthesize a controller based on the bounds of the noise such that a certain target set becomes invariant with respect to the closed-loop dynamics. However, such an approach usually leads to rather conservative controllers and to large infeasibility regions, and although disturbances are not likely to be unbounded in practice, assigning an a priori bound to them seems to demand considerable insight. A stochastic model of the disturbance is a natural alternative approach to this problem: the conservatism of the worst-case analysis may be circumvented, and one need not impose any a priori bounds on the maximum magnitude of the noise. However, since in practice control inputs are almost always bounded, it is of great importance to consider hard bounds on the control inputs as essential ingredients of the controller synthesis; probabilistic constraints on the controllers naturally raise difficult questions on what actions to take when such constraints are violated (see however~\cite{ref:recstrat} for one possible approach to answer these questions).

In this paper we aim to provide answers to the following questions: Given a linear system that is affected by (possibly unbounded) stochastic noise, to be controlled by applying predictive-type {\em bounded} control inputs, (i) is the associated optimization problem tractable? (ii) under what conditions is stability (in a suitable stochastic sense) of the closed-loop system guaranteed? (iii) is stability retained both in the case of MPC implementation and the case of Rolling Horizon Control (RHC) implementation?

In the deterministic setting, there exists a plethora of literature that settles tractability and stability of model-based predictive control, see, for example,  \cite{MayneRawlingsRaoScokaert-00,BemporadMorari-99, ref:maciejowskibk, ref:blanchini1999sic} and the references therein. However, there are fewer results in the stochastic case, some of which we outline next. In~\cite{BertsimasBrown-07}, the authors reformulate the stochastic programming problem as a deterministic one with bounded noise and solve a robust optimization problem over a finite horizon, followed by estimating the performance when the noise can take unbounded values, i.e., when the noise is unbounded, but takes high values with low probability (as in the Gaussian case). In~\cite{Primbs-07,PrimbsSung-08} a slightly different problem is addressed in which the noise enters in a multiplicative manner into the system, and hard constraints on the state and control input are relaxed to probabilistic ones. Similar relaxations of hard constraints to soft probabilistic ones have also appeared in~\cite{ref:CannonKouvaritakisWu-08} for both multiplicative and additive noise inputs, as well as in~\cite{OldewurtelJonesMorari-08}. There are also other approaches, for example those employing randomized algorithms as in~\cite{batinaPhDthesis,MaciejowskiLecchiniLygeros-05}. Finally, a related line of research can be found in \cite{ref:Stoorvogel}, and a novel convex analysis dealing with chance and integrated chance constraints can be found in \cite{AgarwalCinquemaniChatterjeeLygeros-09}.

In this paper we restrict attention to linear time-invariant controlled systems with affine stochastic disturbance inputs. Our approach has three main features. Firstly, for the finite-horizon optimal control subproblem we adopt a feedback control strategy that is affine in certain bounded nonlinear functions of the past noise inputs. Secondly, instead of following the usual trend of adding element-wise constraints to the control input in the optimization, we propose a new approach that entails saturating the utilized noise measurements first and then optimizing over the feedback gains, ensuring that the hard constraints on the input will be satisfied by construction. This novel approach does not require artificially relaxing the hard constraints on the control input to soft probabilistic ones to ensure large feasible sets, and still provides a solution to the problem for a wide class of noise input distributions. In fact, we demonstrate that our strategy (without state constraints) leads to global feasibility. The effect of the noise appears in the finite-horizon optimal control problem as certain covariance matrices, and these matrices may be computed off-line and stored. Thirdly, the measurement saturation functions are only required to be elementwise \emph{bounded} in order to ensure tractability of the optimization problem while maintaining hard constraints on the control input; therefore, these measurement saturation functions may be picked from among the wide class of saturation functions, the standard sigmoidal functions and their piecewise affine approximations, etc.

Once tractability of the finite-horizon underlying optimization problem is insured, it is possible to implement the resulting optimal solution using an MPC approach or an RHC approach. In the former case~\cite{MayneRawlingsRaoScokaert-00}, the optimization problem is resolved at each step and only the first control input is implemented. In the latter case~\cite{AldenSmith-92}, the optimization problem is resolved every $N$ steps (with $N$ being the horizon length) and the entire sequence of $N$ input vectors is implemented. Both of these approaches are shown to provide stability under the assumption that the zero-input and zero-noise system is asymptotically stable, which translates into the condition that the  state matrix $A$ is Schur stable. At a first glance, this assumption might seem restrictive. However, the problem of ensuring bounded variance of linear Gaussian systems with bounded control inputs is, to our knowledge, still open, and here we are considering the problem of controlling a linear system with bounded control input and possibly unbounded noise. It is known that for discrete-time systems without any noise acting on the system it is possible to achieve global stability if and only if the matrix $A$ is neutrally stable~\cite{ref:YangSontagSussmann-97}.

This paper unfolds as follows. In \S\ref{s:ps} we state the main problem to be tackled with the underlying assumptions. In \S\ref{s:if}, we provide a tractable approach to the finite horizon optimization problem with hard constraints on the control input, as well as some examples in \S\ref{sec:examples}. Stability of the MPC and RHC implementations is shown in \S\ref{sec:stability}, and hints onto the input-to-state stable properties of this result are provided in \S\ref{sec:iss}. Finally, we provide a numerical example in \S\ref{sec:Nexample} and conclude in \S\ref{sec:conclusions}.
\subsection*{Notation}
Hereafter, $\NN := \{1,2,\ldots\}$ is the set of natural numbers, $\Nz := \NN \cup \{0\}$, and $\posR$ is the set of nonnegative real numbers. We let $\indic{A}(\cdot)$ denote the indicator function of a set $A$, and $\mathbf I_{n\times n}$ and $\mathbf 0_{n\times n}$ denote the $n$-dimensional identity and zeros matrices, respectively. Also, let $\EE_{\xz}[\cdot]$ denote the expected value given $\xz$, and $\tr{\cdot}$ denote the trace of a matrix. For a given symmetric $n$-dimensional matrix $M$ with real entries, let $\{\lambda_i(M)\mid i=1, \ldots, n\}$ be the set of eigenvalues of $M$, and let $\lambda_{\rm max}(M) := \max_i\lambda_i(M)$ and $\lambda_{\text{min}}(M) := \min_i\lambda_i(M)$. Let  $\norm{\cdot}_p$ denote standard $\ell_p$ norm. Finally, the mean and covariance matrix of any vector $v$ are denoted by $\Sigma_v$ and $\mu_v$, respectively.

\section{Problem Statement}
\label{s:ps}
Consider the following general affine discrete-time stochastic dynamical model:
\begin{equation}
\label{eq:system}
	x_{t+1} = A x_t + Bu_t + Fw_t+r, \qquad t\in\NN_0,
\end{equation}
where $x_t\in\R^n$ is the state, $u_t\in \R^m$ is the
control input, $w_t\in\R^n$ is a stochastic noise input vector, $A$, $B$ and $F$ are known matrices, and $r\in\R^n$ is a known constant vector. We assume that the  initial condition $x_0$ is given and
that, at any time $t$, $x_t$ is observed exactly. We shall assume further that the noise vectors $w_t$ are i.i.d. and that the control input vector is bounded at each instant of time $t$, i.e.,
\begin{equation}
\label{eq:bddu}
	u_t\in\mathbb{U} := \bigl\{u\in\mathbb R^m\big| \norm{u}_\infty \leq U_{\rm max}\bigr\} \quad \fa t\in\Nz,
\end{equation}
where $U_{\mathrm{max}}>0$ is some given element-wise saturation bound. Note that the model \eqref{eq:system} with constraints \eqref{eq:bddu} can handle a wide range of convex polytopic constraints. In particular, any system
\begin{equation}\label{eqn:generalsystem}
    x_{t+1}=Ax_t+\hat Bv_t+F\hat w_t+\hat r
\end{equation}
with input constraints $v_t\in\mathbb V$ that can be transformed to the form \eqref{eq:bddu} by an affine transformation \[v_t=Su_t+l\]
is amenable to our approach by setting $B=\hat BS$ and $r=\hat Bl+\hat r$ in \eqref{eq:system}. Note that the set $\mathbb V$ need not necessarily be a hypercube, or even contain the origin. Note also that  we can assume that $w_t$ is zero mean in \eqref{eq:system} without loss of generality; given a system of the form \eqref{eqn:generalsystem} where $\hat w_t$ is not zero mean, we can  replace it by a system in the form \eqref{eq:system} with zero mean in which
\[w_t=\hat w_t -\EE[w_t]\]
by setting $r=\hat r+F\EE[w_t]$.


Fix a horizon $N\in\NN$ and set $t = 0$. The \emph{MPC} procedure can be described as follows.
\begin{itemize}
	\item[(a)] Determine an admissible optimal feedback control policy, say $\pi^\star_{t:t+N-1}\in\Pi$, for an $N$-stage cost function starting from time $t$, given the (measured) initial condition $x_t$; 
	\item[(b)] increase $t$ to $t+1$, and go back to step (a).
\end{itemize}
On the other hand, the \emph{RHC} procedure simply replaces (b) above by
\begin{itemize}
	\item[(b$'$)] apply the entire sequence $\pi^\star_{t:t+N-1}$ of control inputs, update the state $x_{t+N}$ at the $(t+N-1)$-th step, increase $t$ to $t+N$ and go back to step (a).
\end{itemize}
Accordingly, the $t$-th step of this procedure consists of minimizing the stopped $N$-period cost function starting at time $t$, namely, the objective is to find a feedback control policy that attains
\begin{align}
	\inf_{\pi\in\Pi} V_{t, t+N-1}(\pi, x) :=& \;\inf_{\pi\in\Pi}\mathsf E^\pi_{x_t}\!\biggl[\sum_{i=t}^{t+N-1} \!\!\bigl(x_i\transp Q_i x_i+u_i\transp R_i u_i\bigr) \nn\\ &\qquad\qquad+ x_{t+N}\transp Q_{t+N}x_{t+N}\biggr].\label{e:rhcf}
\end{align}
Since both the system \eqref{eq:system} and cost \eqref{e:rhcf} are time-invariant, it is enough to consider the problem of minimizing the cost for $t = 0$, i.e., the problem of minimizing $V_{0, N-1}(\pi, x)$ over $\pi\in\Pi$.\\
In view of the above we consider the problem
\begin{equation}
\label{eq:problem}
\begin{aligned}
\min_{\pi\in\Pi}   \quad& \EE_{\xz}\left[\sum\limits_{t=0}^{N-1}\bigl(x_t\transp Q_t x_t+u_t\transp R_t u_t\bigr)+x_N\transp Q_Nx_N\right], \\
\textrm{s.t.}\quad&\mathrm{dynamics}\,(\ref{eq:system}),\, \mathrm{and\,\, constraints}\,(\ref{eq:bddu})
\end{aligned}
\end{equation}
where $Q_t>0$ and $R_t>0$ are some given symmetric matrices of appropriate dimension. If feasible with respect to \eqref{eq:bddu}, Problem \eqref{eq:problem} generates an optimal sequence of feedback control laws $\pi^*=\{u^*_0,\cdots, u^*_{N-1}\}$.

The evolution of the system (\ref{eq:system}) over a single optimization horizon $N$ can be described in compact form as follows:
\begin{equation}
\label{eq:compactdyn}
\bar{x}=\bar{A}\x0+\bar{B}\bar{u}+\bar{D}\bar F\bar{w}+\bar D\bar r,
\end{equation}
where
$$\bar{x}:= \begin{bmatrix}
x_0 \\ x_1  \\ \vdots \\ x_N
\end{bmatrix},\, \bar{u}:= \begin{bmatrix}
u_0 \\ u_1 \\ \vdots \\ u_{N-1}
\end{bmatrix},\, \bar r :=\smat{r\\ \vdots \\ r}, \,
\bar{w} :=
\begin{bmatrix}
w_0 \\ w_1 \\ \vdots \\ w_{N-1}
\end{bmatrix},$$
$$
\bar{A} :=
\begin{bmatrix}
\mathbf I_{n\times n} \\ A \\ \vdots \\ A^N
\end{bmatrix},\, \bar{B} :=
\begin{bmatrix}
\mathbf 0_{n\times m} &\cdots &\cdots & \mathbf 0_{n\times m} \\
B &\ddots &&\vdots\\
AB & B &\ddots & \vdots\\
\vdots && \ddots &\mathbf 0_{n\times m}\\
A^{N-1} B & \cdots & AB& B
\end{bmatrix},$$ $$
{\small \bar{D} :=
\begin{bmatrix}
\mathbf 0_{n\times n} &  \cdots &\cdots& \mathbf 0_{n\times n} \\
\mathbf I_{n\times n} & \ddots & &\vdots \\
A & \mathbf{I}_{n\times n} &\ddots &\vdots \\
\vdots && \ddots & \mathbf 0_{n\times n} \\
A^{N-1} &\cdots & A& \mathbf I_{n\times n}
\end{bmatrix}, \, \bar F:=\smat{F&\hdots & \mathbf 0 \\ \vdots &\ddots &\vdots \\  \mathbf 0 &\hdots & F}}
$$
where the input
\begin{equation}\label{eq:bddu2}
\bar u\in\bar{\mathbb{U}}:= \bigl\{\xi\in\mathbb R^{Nm}\big| \norm{\xi}_\infty \leq U_{\rm max}\bigr\}.
\end{equation}
Using the compact notation above, the optimization Problem~\eqref{eq:problem} can be rewritten as follows:
\begin{equation}
\label{eq:problem1}
\begin{aligned}
\min_{\pi\in\Pi}   \quad&\EE_{\xz}\bigl[\bar x\transp \bar Q\bar x+\bar u\transp \bar R\bar u\bigr], \\
\textrm{s.t.}\quad& {\rm dynamics}\, (\ref{eq:compactdyn}),\, \mathrm{and\,\, constraints}\,(\ref{eq:bddu2}) ,\\
\end{aligned}
\end{equation}
where
\[\bar Q=\smat{Q_0&\hdots & \mathbf 0_{n\times n} \\ \vdots& \ddots& \vdots\\  \mathbf 0_{n\times n} &\hdots& Q_N },\,
\bar R=\smat{R_0&\hdots & \mathbf 0_{m\times m} \\
\vdots& \ddots& \vdots\\  \mathbf 0_{m\times m} &\hdots& R_{N-1} }.\]

The solution to Problem \eqref{eq:problem1} is difficult to obtain in general.  In order to obtain an optimal solution to Problem \eqref{eq:problem1} over the class of feedback policies, we need to solve the Dynamic Programming equations. This generally requires using some gridding technique, making the problem  extremely difficult to solve computationally. Another approach is to restrict attention to a specific class of state feedback policies. This will  result in a suboptimal solution to our problem, but may yield a tractable optimization problem. It is the track we pursue in the next section.

\section{Tractable Solution under Bounded Control Inputs}
\label{s:if}
By the hypothesis that the state is observed without error, one
may reconstruct the noise sequence from the sequence of observed
states and inputs by the formula
\begin{equation}
\label{eq:noiserec}
Fw_t=x_{t+1}-A x_t -B u_t-r,\qquad t\in\NN_0.
\end{equation}
In the light of this, and inspired by the works~\cite{ref:ben-tal04, ref:goulart06}, we
shall consider feedback policies of the form:
\begin{equation}
\label{eq:controlpolicy}
u_t=\sum_{i=0}^{t-1} G_{t,i} F w_i + d_t,
\end{equation}
where the feedback gains $G_{t,i}\in\R^{m\times n}$ and the
affine terms $d_t\in\R^m$ must be chosen based on the control
objective, while observing the constraints \eqref{eq:bddu}. With this definition, the value of $u$ at time $t$
depends on the values of $w$ up to time $t-1$.
Using~(\ref{eq:noiserec}) we see that $u_t$ is a function of the
observed states up to time $t$. It was shown in~\cite{ref:goulart06}
that there exists a one-to-one (nonlinear) mapping between control
policies in the form~(\ref{eq:controlpolicy}) and the class of
affine state feedback policies. That is, provided one is
interested in affine state feedback policies,
parametrization~(\ref{eq:noiserec}) constitutes no loss of
generality. Of course, this choice is generally suboptimal, but it
will ensure the tractability of a large class of optimal control
problems. In compact notation, the control sequence up to time $N-1$
is given by
\begin{equation}\label{eq:auginput}
\bar{u}=\bar{G}\bar F\bar{w}+\bar{d},
\end{equation}
where $\bar{d} :=
\smat{
d_0\transp & d_2\transp & \hdots & d\transp_{N-1}}\transp$, and
\[\bar{G} :=
\begin{bmatrix}
\mathbf 0_{m\times n} \\
G_{1,0} &\mathbf 0_{m\times n} \\
\vdots &\ddots &\ddots \\
G_{N-1,0} &\cdots &  G_{N-1,N-2} &\mathbf 0_{m\times n}
\end{bmatrix}. \]

Since the elements of the noise vector $\bar w$ are not assumed to be bounded, there can be no guarantee that the control input \eqref{eq:auginput} will meet the constraint \eqref{eq:bddu2}. This is a problem in practical applications, and has traditionally been circumvented by assuming that the noise input lies within a compact set~\cite{ref:goulart06}, and designing a worst-case controller. In this article we propose to use the controller
\begin{equation}\label{eq:auginputbdd}
\bar{u}=\bar{G}\bar\varphi(\bar F\bar{w})+\bar{d},
\end{equation}
instead of (\ref{eq:auginput}), where
\[
	\bar \varphi(\bar F\bar w)=\smat{\varphi_0(F{w}_0)\\ \vdots \\ \varphi_{N-1}(F{w}_{N-1})},
\]
$\varphi_i(Fw_i)$ is a shorthand for the vector $\bigl[\varphi_i^1(F_1w_i), \ldots, \varphi_i^n(F_nw_i)\bigr]\transp$, $F_j$ is the $j$-th row of the matrix $F$, and $\varphi_i^j:\R\to\R$ is any function with $\sup\limits_{s\in \mathbb{R}}|\varphi_i^j(s)|\leq \phi_{\max} \leq U_{\rm max}$. In other words, we have chosen to saturate the measurements that we obtain from the noise input vector before inserting them into our control vector. This way we do not assume that the noise distribution is defined over a compact domain, which is an advantage over other approaches \cite{BertsimasBrown-07,ref:goulart06}.
Moreover, the choice of element-wise saturation functions $\varphi_i(\cdot)$ is left open. As such, we can accommodate standard saturation, piecewise linear, and sigmoidal functions, to name a few.

\begin{remark}
	Our choice of saturating the measurement from the noise vectors renders the optimization problem tractable as opposed to just calculating the whole input vector $\bar u$ and then saturating it afterwards, which tends to an intractable optimization problem.\RemEnd
\end{remark}

\begin{remark}
	Note that the choices of control inputs in \eqref{eq:auginput} and \eqref{eq:auginputbdd} are both \emph{non Markovian}; however, they differ in the fact that the former depends affinely on previous noise inputs $\bar w$, whereas the latter is a nonlinear feedback due to passing noise measurements through the function $\bar \varphi(.)$.\RemEnd
\end{remark}

\begin{proposition}
\label{p:main}
Assume that $\EE_{\xz}\left[\bar\varphi(\bar F\bar w)\right]=0$, $\forall x_0\in\R^n$. Then, Problem~\eqref{eq:problem1} with the input~\eqref{eq:auginputbdd} is a convex optimization problem, with
respect to the decision variables $(\bar G,\bar d)$, which is given by
\begin{equation}
\label{eq:problem2}
\begin{aligned}
\min\limits_{(\bar G,\bar d)} \quad& b\transp \bar d + \bar d\transp M_1\bar d + \tr{\bar G\transp M_1\bar G\Lambda_1 + M_2\bar G\Lambda_2} \\
\mathrm{s.t.}\quad& |\bar d_i|+ \norm{\bar G_i}_1\phi_{\rm max}  \leq U_{\max},\quad \forall i=1,\cdots,Nm
\end{aligned}
\end{equation}
where $G_i$ is the $i$-th row of $G$,
\begin{align*}
  b^T &= 2(\bar A x_0+\bar D\bar F\mu_{\bar w}+\bar r)\transp \bar Q\bar B, \quad
  M_1 = \bar R+ \bar B\transp \bar Q\bar B, \\
  M_2 &= 2\bar F\transp\bar D\transp\bar Q\bar B, \\
  \Lambda_1 &= \mathrm{diag}\bigl\{\EE\bigl[\varphi_0(Fw_0)\varphi_0(Fw_0)\transp\bigr],\cdots, \bigr. \\ &\qquad \qquad \bigl. \EE\bigl[\varphi_{N-1}(Fw_{N-1})\varphi_{N-1}(Fw_{N-1})\transp\bigr]\bigr\},\\
  \Lambda_2 &= \mathrm{diag}\bigl\{\EE\bigl[\varphi_0(Fw_0) w_0\transp\bigr],\cdots,\bigr. \\ &\qquad \qquad \bigl.\EE\bigl[\varphi_{N-1}(Fw_{N-1})w_{N-1}\transp\bigr]\bigr\}.
\end{align*}
\end{proposition}
\begin{proof}
Let us first consider the cost function in Problem \ref{eq:problem1}.
After substituting the system equations, we obtain
\begin{align}
  &\EE_{\xz}\bigl[\bar{x}\transp \bar Q \bar x +\bar{u}\transp \bar R\bar u\bigr] =\\ & \EE_{\xz}[\left(\bar A x_0 +\bar B \bar u+\bar D\bar F\bar w+\bar r\right)\transp \bar Q\left(\bar A x_0 +\bar B \bar u+\bar D\bar F\bar w+\bar r\right)\nn\\ &\qquad +\bar u\transp \bar R\bar u] \nn\\
  &=(\bar A x_0+\bar r)\transp \bar Q(\bar A x_0+\bar r) + 2(\bar A x_0+\bar r)\transp \bar Q\bar D \bar F\EE_{\xz}\bigl[\bar w\bigr]\nn\\ &\quad  + 2(\bar A x_0+\bar r)\transp\bar Q\bar B \EE_{\xz}\bigl[\bar u\bigr]   +2\EE_{\xz}\bigl[ \bar w\transp \bar F\transp\bar D\transp \bar Q \bar B\bar u\bigr] \nn\\
  &\quad +\EE_{\xz}\bigl[ \bar w \transp \bar F\transp\bar D\transp \bar Q\bar D \bar F\bar w\bigr] +\EE_{\xz}\bigl[ \bar u\transp (\bar R+ \bar B\transp \bar Q \bar B)\bar u\bigr].\nn
\end{align}
Note that since $\EE_{\xz}\left[\bar\varphi(\bar F\bar w)\right]=0$, we have that $\EE_{\xz}\bigl[\bar u\bigr]=\bar d$. Accordingly, using the definitions of $b$, $M_1$, $M_2$, and $\Lambda_2$,
\begin{align}
  \EE_{\xz}\bigl[\bar{x}\transp \bar Q \bar x +\bar{u}\transp \bar R\bar u\bigr]
  &= b\transp \bar d+\tr{M_2\bar G\Lambda_2}  +c \nn\\ &\quad + \EE_{\xz}\bigl[ \bar u\transp M_1\bar u\bigr],\label{eq:arg0}
\end{align}
where $  c= (\bar A x_0+\bar r)\transp \bar Q(\bar A x_0+\bar r)
   + \tr{\bar F\transp\bar D\transp \bar Q\bar D \bar F\Sigma_{\bar w}} +2(\bar A x_0+\bar r)\transp \bar Q\bar D\bar F\mu_{\bar w}$ is a constant that we omit as it does not change the optimization problem, and we have used the following intermediate step
\begin{align*}
   &\EE_{\xz}\bigl[ \bar w\transp \bar F\transp\bar D\transp \bar Q \bar B\bar u\bigr]
   = \EE_{\xz}\bigl[ \bar w\transp \bar F\transp\bar D\transp \bar Q \bar B(\bar G\bar \varphi(\bar F\bar w)+\bar d)\bigr]\\ &\qquad =  \tr{\bar F\transp\bar D\transp \bar Q\bar B\bar G\Lambda_2}+\mu_{\bar w}\transp \bar F\transp\bar D\transp \bar Q \bar B\bar d.
\end{align*}
Using again the assumption that $\EE_{\xz}\left[\bar\varphi(\bar F\bar w)\right]=0$, we have that
\begin{align}
  \EE_{\xz}\bigl[ \bar u\transp & M_1\bar u\bigr] = \EE_{\xz}\bigl[(\bar{G}\bar\varphi(\bar F\bar{w})+\bar{d})\transp M_1(\bar{G}\bar\varphi(\bar F\bar w)+\bar{d})\bigr] \nn\\
   &= \EE_{\xz}\bigl[\bar\varphi(\bar F\bar{w})\transp \bar G\transp M_1\bar G\bar\varphi(\bar F\bar{w})\bigr]+\bar d\transp M_1\bar d \nn\\
   &= \tr{\bar G\transp M_1\bar G \EE_{\xz}\bigl[\bar\varphi(\bar F\bar{w})\bar\varphi(\bar F\bar{w})\transp \bigr]}+\bar d\transp M_1\bar d \nn \\
   &= \tr{\bar G\transp M_1\bar G\Lambda_1}+\bar d\transp M_1\bar d. \label{eq:arg2}
\end{align}
Finally, combining \eqref{eq:arg0} and \eqref{eq:arg2}, we obtain the cost in Problem \ref{eq:problem2}, which is convex.

Let us look at the constraints in Problem \ref{eq:problem1}. The proposed control input  \eqref{eq:auginputbdd} satisfies the hard constraints \eqref{eq:bddu2} as long as the following condition is satisfied: $\norm{\bar d+\bar G\bar\varphi(\bar w)}_\infty\leq U_{\max}$,  $\forall \bar\varphi(\bar w)$ such that $\norm{\bar\varphi(\bar w)}_\infty\leq \phi_{\max}$. This is equivalent to the following  conditions: $\forall i=1,\cdots, Nm$, $|\bar d_i+\bar G_i\bar\varphi(\bar w)|\leq U_{\max}$, $\forall \bar\varphi(\bar w)$ such that $\norm{\bar\varphi(\bar w)}_\infty\leq \phi_{\max}$. As these conditions should hold for any permissible value of the function $\bar\varphi(\bar w)$, we can eliminate the dependence of the constraints on $\bar\varphi(\bar w)$ through the following optimization problems
$\max\limits_{\norm{\bar\varphi(\bar w)}_\infty\leq \phi_{\max}}|\bar d_i+\bar G_i\bar\varphi(\bar w)|\leq U_{\max},\, \forall i=1,\cdots, Nm$.
It is straightforward now to show, using H\" older's inequality \cite[p.~29]{Luenberger-69}, that $\max\limits_{\norm{\bar\varphi(\bar w)}_\infty\leq \phi_{\max}}|\bar d_i+\bar G_i\bar\varphi(\bar w)| = |\bar d_i|+ \norm{\bar G_i}_1\phi_{\rm max}$,
and the result follows.
\end{proof}

\begin{remark}
	Problem~\eqref{eq:problem2} is a quadratic program in the optimization parameters $\theta := (\bar G, \bar d)$~\cite[p.~111]{ref:boyd04}, and can be solved efficiently by standard solvers such as \texttt{cvx}~\cite{ref:boydCVX}.\RemEnd
\end{remark}

\subsection{Examples}\label{sec:examples}
	An important step in the solvability of Problem \eqref{eq:problem2} is being able to calculate the matrices $\Lambda_1$ and $\Lambda_2$. In general, these matrices can be calculated off-line by numerical integration. However, in some instances these matrices can be given in terms of explicit formulas; two of these instances are given in the following examples.

Recall the following standard special mathematical functions: the \emph{standard error function} $\erf(z) := \frac{2}{\sqrt\pi}\int_0^z \epower{-\frac{t^2}{2}}\mrm dt$ and the \emph{complementary error function}~\cite[p.~297]{ref:AbramowitzStegun} defined by $\erfc(z) := 1-\erf(z)$ for $z\in\R$, the \emph{incomplete Gamma function}~\cite[p.~260]{ref:AbramowitzStegun} defined by $\Gamma(a, z) := \int_z^\infty t^{a-1} \epower{-t} \mrm dt$ for $z, a > 0$, the \emph{confluent hypergeometric function}~\cite[p.~505]{ref:AbramowitzStegun} defined by $U(a, b, z) := \frac{1}{\Gamma(a)}\int_0^\infty \epower{-zt} t^{a-1} (1+t)^{b-a-1}\mrm dt$ for $a, b, z > 0$ and $\Gamma$ is the standard Gamma function.

	We collect a few facts in the following
	\begin{proposition}
	\label{p:collect}
		For $\sigma^2 > 0$ we have
		\begin{enumerate}
			\item $\displaystyle{\frac{1}{\sqrt{2\pi}\sigma}\int_z^\infty \epower{-\frac{t^2}{2\sigma^2}}\mrm dt = \frac{1}{2}\Bigl(1 + \erf\Bigl(\frac{z}{\sqrt 2\sigma}\Bigr)\Bigr)}$;
			\item $\displaystyle{\frac{1}{\sqrt{2\pi}\sigma}\int_0^\infty \frac{t^2}{1+t^2} \epower{-\frac{t^2}{2\sigma^2}}\mrm dt}$ \\$\displaystyle{\qquad =\frac{1}{2}\Bigl(\sqrt{2\pi}\sigma - \pi\epower{-\frac{1}{2\sigma^2}}\erfc\Bigl(\frac{1}{\sqrt 2\sigma}\Bigr)\Bigr)}$;
			\item $\displaystyle{\frac{1}{\sqrt{2\pi}\sigma}\int_0^1 t^2 \epower{-\frac{t^2}{2\sigma^2}}\mrm dt}$\\ $\displaystyle{\qquad= \sqrt{\frac{\pi}{2}}\sigma^3\erf\Bigl(\frac{1}{\sqrt 2\sigma}\Bigr) - \sigma^2\epower{-\frac{1}{2\sigma^2}}}$;
			\item $\displaystyle{\frac{1}{\sqrt{2\pi}\sigma}\int_1^\infty t \epower{-\frac{t^2}{2\sigma^2}}\mrm dt = \frac{\sigma}{\sqrt{2\pi}}\Gamma(2\sigma^2, 1)}$;
			\item $\displaystyle{\frac{1}{\sqrt{2\pi}\sigma}\int_0^\infty \frac{t^2}{\sqrt{1+t^2}} \epower{-\frac{t^2}{2\sigma^2}}\mrm dt = \frac{\sigma}{2\sqrt 2}U\Bigl(\frac{1}{2}, 0, \frac{1}{2\sigma^2}\Bigr)}$.
		\end{enumerate}
	\end{proposition}

	\begin{example}
	\label{ex:sigmoids}
		Let us consider~\eqref{eq:system} with Gaussian noise and sigmoidal bounds on the control input. More precisely, suppose that the noise process $(w_t)_{t\in\Nz}$ is an independent and identically distributed (i.i.d) sequence of Gaussian random vectors of mean $0$ and covariance $\Sigma$. Let the components of $w_t$ be mutually independent, which implies that $\Sigma$ is a diagonal matrix $\diag\{\sigma_1^2, \ldots, \sigma_n^2\}$. Suppose further that the matrix $F=I$ and that the function $\varphi$ is a standard sigmoid, i.e., $\varphi(t) := t/\sqrt{1+t^2}$. Then from Proposition~\ref{p:collect} we have for $i=1, \ldots, n$ and $j=0, \ldots, N-1$,
		\begin{align*}
			\EE[\varphi(w_j^i)^2] 
			&= 2\cdot\frac{1}{\sqrt{2\pi}\sigma_i}\int_{0}^\infty \frac{t^2}{1+t^2} \epower{-\frac{t^2}{2\sigma_i^2}}\\
			& = \sqrt{2\pi}\sigma_i - \pi\epower{-\frac{1}{2\sigma_i^2}}\erfc\Bigl(\frac{1}{\sqrt 2\sigma_i}\Bigr).
		\end{align*}
		This shows that the matrix $\Lambda_1$ in Proposition~\ref{p:main} is equal to $\diag\{\Sigma', \ldots, \Sigma'\}$, where $
			\Sigma' := \diag\left\{\sqrt{2\pi}\sigma_1 - \pi\epower{-\frac{1}{2\sigma_1^2}}\erfc\Bigl(\frac{1}{\sqrt 2\sigma_1}\Bigr),\right.$ $ \left. \ldots, \sqrt{2\pi}\sigma_n - \pi\epower{-\frac{1}{2\sigma_n^2}}\erfc\Bigl(\frac{1}{\sqrt 2\sigma_n}\Bigr)\right\}.
$
		Similarly, since
		\begin{align*}
			\EE[\varphi(w_j^i) w_j^i] 
			&=  \frac{2}{\sqrt{2\pi}\sigma_i}\int_{-\infty}^\infty \frac{t^2}{\sqrt{1+t^2}} \epower{-\frac{t^2}{2\sigma_i}} \mrm dt\\
			& = \frac{\sigma_i}{\sqrt 2}U\Bigl(\frac{1}{2}, 0, \frac{1}{2\sigma_i^2}\Bigr),
		\end{align*}
		the matrix $\Lambda_2$ in Proposition~\ref{p:main} is $\diag\{\Sigma'', \ldots, \Sigma''\}$, where
		$
			\Sigma'' := \diag\left\{\frac{\sigma_1}{\sqrt 2}U\Bigl(\frac{1}{2}, 0, \frac{1}{2\sigma_1^2}\Bigr), \ldots, \frac{\sigma_n}{\sqrt 2}U\Bigl(\frac{1}{2}, 0, \frac{1}{2\sigma_n^2}\Bigr)\right\}.
		$
		Therefore, given the system~\eqref{eq:system}, the control policy~\eqref{eq:controlpolicy}, and the description of the noise input as above, the matrices $\Lambda_1$ and $\Lambda_2$ derived above complete the set of hypotheses of Proposition~\ref{p:main}. The problem~\eqref{eq:problem} can now be solved as a quadratic program~\eqref{eq:problem2}.\ExEnd
	\end{example}

Note that we have chosen to use the standard sigmoidal functions in Example \ref{ex:sigmoids}. However, the result still holds for more general sigmoidal functions of the form $\tilde\phi(t)=M\frac{\alpha t}{\sqrt{1+\alpha^2t^2}}$, where $M\in\R$ is some given magnitude and $\alpha\in \R$ is some given slope. This slight change is reflected in the entries of the matrices $\Lambda_1$ and $\Lambda_2$, i.e.,
for $i=1, \ldots, n$ and $j=0, \ldots, N-1$,
		\begin{align*}
			\EE[\varphi(w_j^i)^2] & = M\left(\sqrt{2\pi}\sigma_i\alpha - \pi\epower{-\frac{1}{2\sigma_i^2\alpha^2}}\erfc\left(\frac{1}{\sqrt 2\sigma_i\alpha}\right)\right),
		\end{align*}
and	$\EE[\varphi(w_j^i) w_j^i]  = M\frac{\sigma_i \alpha}{\sqrt 2}U\Bigl(\frac{1}{2}, 0, \frac{1}{2\sigma_i^2\alpha^2}\Bigr)$.

	\begin{example}\label{ex:saturation}
		Consider the system~\eqref{eq:system} as in Example~\ref{ex:sigmoids}, with $\varphi$ being the standard saturation function defined as $\varphi(t) = \sat(t) := \sgn(t)\min\{|t|, 1\}$. From Proposition~\ref{p:main} we have for $i=1, \ldots, n$ and $j=0, \ldots, N-1$,
		\begin{align*}
			\xi_i' &:= \EE[\varphi(w_j^i)^2] = \frac{1}{\sqrt{2\pi}\sigma_i}\int_{-\infty}^\infty \varphi(t)^2 \epower{-\frac{t^2}{2\sigma_i^2}} \mrm dt\\
			& = \frac{2}{\sqrt{2\pi}\sigma_i}\int_0^1 t^2\epower{-\frac{t^2}{2\sigma_i^2}} \mrm dt + \frac{2}{\sqrt{2\pi}\sigma_i}\int_1^\infty \epower{-\frac{t^2}{2\sigma_i^2}} \mrm dt\\
			& = \sqrt{2\pi}\sigma_i^3\erf\Bigl(\frac{1}{\sqrt 2\sigma_i}\Bigr) - 2\sigma_i^2\epower{-\frac{1}{2\sigma_i^2}} + 1 + \erf\Bigl(\frac{1}{\sqrt 2\sigma_i}\Bigr)
		\end{align*}
		and
		\begin{align*}
			\xi_i'' & := \EE[\varphi(w_j^i) w_j^i] = \frac{1}{\sqrt{2\pi}\sigma_i}\int_{-\infty}^\infty t \varphi(t) \epower{-\frac{t^2}{2\sigma_i^2}} \mrm dt\\
			& = \frac{2}{\sqrt{2\pi}\sigma_i}\int_0^1 t^2\epower{-\frac{t^2}{2\sigma_i^2}} \mrm dt + \frac{2}{\sqrt{2\pi}\sigma_i}\int_1^\infty t\epower{-\frac{t^2}{2\sigma_i^2}} \mrm dt\\
			& = \sqrt{2\pi}\sigma_i^3\erf\Bigl(\frac{1}{\sqrt 2\sigma_i}\Bigr) - 2\sigma_i^2\epower{-\frac{1}{2\sigma_i^2}} + \sqrt{\frac{2}{\pi}}\sigma_i\Gamma(2\sigma_i^2, 1).
		\end{align*}
		Therefore, in this case the matrix $\Lambda_1$ in Proposition~\ref{p:main} is $\diag\{\Sigma', \ldots, \Sigma'\}$ with $\Sigma' := \diag\{\xi_1', \ldots, \xi_n'\}$, and the matrix $\Lambda_2$ is $\diag\{\Sigma'', \ldots, \Sigma''\}$ with $\Sigma'' := \diag\{\xi_1'', \ldots, \xi_n''\}$. These information complete the set of hypotheses of Proposition~\ref{p:main}, and the problem~\eqref{eq:problem} can now be solved as a quadratic program~\eqref{eq:problem2}.\ExEnd
	\end{example}
\section{Stability Analysis} \label{sec:stability}
In this section, we assume that the matrix $A$ is Schur stable, i.e., $\abs{\lambda_i(A)}<1$, $\fa i$. Accordingly, and since the control is bounded, it is intuitively evident that the closed-loop system is stable in some sense. Indeed, we shall show that the variance of the state is uniformly bounded both in the MPC and RHC cases, the only difference being a choice of implementation based on available memory.

First we need the following Lemma. It is a standard variant of the Foster-Lyapunov condition~\cite{ref:meyn93}; we include a proof here for completeness. The hypotheses of this Lemma are stronger than usual, but are sufficient for our purposes; see e.g.,~\cite{ref:foss04} for more general conditions.

	\begin{lemma}
	\label{l:FL}
		Let $(x_t)_{t\in\Nz}$ be an $\R^n$-valued Markov process. Let $V:\R^n\to\posR$ be a continuous positive definite and radially unbounded function, integrable with respect to the probability distribution function of $w$. Suppose that there exists a compact set $K \subset\R^n$ and a number $\lambda \in\;]0, 1[$ such that
			\[
				\EE\bigl[V(x_{1}) \big| x_0 = x\bigr] \le \lambda V(x),\qquad \forall x\not\in K.
			\]
			Then $\sup\limits_{t\in\Nz} \EE_{x}\bigl[V(x_t)\bigr] < \infty$.
	\end{lemma}
	\begin{proof}
From the conditions it follows immediately that
		\[
			\EE_{x}\bigl[V(x_{1})\bigr] \le \lambda V(x) + b\indic{K}(x), \qquad \fa x\in \mathbb R^n
		\]
		where $b := \sup\limits_{x\in K}\EE_x\bigl[V(x_{1})\bigr]$. We then have
		\begin{align}
		\label{eq:Vtbound}
			\EE_x\bigl[V(x_{t})\bigr] & = \EE_x\bigl[\EE\bigl[V(x_{t})\big|x_{t-1}\bigr]\bigr] \\ & \le \EE_x\bigl[\EE\bigl[\lambda V(x_{t-1}) + b\indic{K}(x_{t-1})\bigr]\bigr]\nn\\
			& \le \lambda^t V(x) + \sum_{i=0}^{t-1}\lambda^{t-1-i} b\;\EE_x\bigl[\indic{K}(x_i)\bigr]\nn\\
			& \le \lambda^t V(x) + \frac{b(1-\lambda^{t})}{1-\lambda},
		\end{align}
		which shows that $\sup\limits_{t\in\Nz}\EE_{x}\bigl[V(x_t)\bigr] \le V(x) + b/(1-\lambda) < \infty$ as claimed.
	\end{proof}

We shall utilize Lemma \ref{l:FL} in order to show that the implementation of either the MPC or the RHC  strategy generated by the solution of Problem \eqref{eq:problem2} results in a uniformly bounded state variance. 
\subsection{MPC Case}
The MPC implementation corresponding to our input~\eqref{eq:auginputbdd} and optimization program~\eqref{eq:problem2} consists of the following steps: Given a fixed optimization horizon $N$, set the initial time $t=0$,  calculate the optimal control gains $(\bar G^*,\bar d^*)$ using the program \eqref{eq:problem2}, apply the first optimal control input $\pi_{0|t}^*=u_{0|t}^*=\bar d_{0|t}^*$, increase $t$ to $t+1$, and iterate. Of course, the optimal gain depends implicity on the current given initial state, i.e., $\bar d_{0|t}^*=\bar d_{0|t}^*(x_t)$, which in turn gives rise to a stationary infinite horizon optimal policy given by $\mathbf\pi^{\rm MPC} := \bigl(\pi_{0|0}^*, \pi_{0|1}^*, \ldots\bigr)=\bigl(\bar d_{0|t}^*,\bar d_{0|t}^*, \ldots\bigr)$.
 The closed-loop system is thus given by
	\begin{equation}
	\label{eq:clsys}
		x_{t+1} = Ax_t + B \bar d^*_{0|t} + Fw_t+r,\qquad t\in\Nz.
	\end{equation}
	\begin{proposition}
	\label{prop:SHmain}
		Assume that the matrix $A$ is Schur stable and the assumptions of Proposition \ref{p:main} hold. Then, under the control policy $\pi^{\rm MPC}$ defined above, the closed loop system~\eqref{eq:clsys} satisfies $\sup_{t\in\Nz} \EE_{\xz}\Bigl[\norm{x_t}^2\Bigr] < \infty$.
	\end{proposition}
	\begin{proof}
		Since by assumption the matrix $A$ is Schur stable, there exists a positive definite and symmetric matrix with real entries, say $P$, such that $A\transp P A - P \le -\mathbf{I}_{n\times n}$. Using the system \eqref{eq:clsys}, at each time instant $t\in \Nz$ we have
		\begin{align*}
			&\EE_{x_t}\bigl[x_{t+1}\transp Px_{t+1}\bigr] = \\ & \EE_{x_t}\bigl[(Ax_t + B \bar d^*_{0|t} + Fw_t+r)\transp P(Ax_t + B \bar d^*_{0|t} + Fw_t+r)\bigr]\\
				& = x_t\transp A\transp P A x_t + 2 x_t\transp A\transp P(B\bar d^*_{0|t}+F\mu_{w_t}+r) \\ & \quad + \bar d^{*\mathsf T}_{0|t} B\transp PB\bar d^*_{0|t} +r\transp P r + 2(F\mu_{w_t}+r)\transp PB \bar d^*_{0|t}\\ &\quad  + 2r\transp PF\mu_{w_t}
+ \tr{F\transp PF\Sigma_{w_t}}.
		\end{align*}
			Using the fact that $\norm{\bar d^*_{0|t}}_\infty \le U_{\text{max}}$ (from ~\eqref{eq:problem2}), we obtain the following bound
        \begin{align*}
         \EE_{x_t}\bigl[x_{t+1}\transp Px_{t+1}\bigr] & \leq  x_t\transp A\transp P A x_t +2c_1\norm{x_t}_\infty +c_2,
        \end{align*}
where $c_1:=\norm{A\transp P(F\mu_{w_t}+r)}_1+m\norm{A\transp PB}_\infty U_{\max}$ and $c_2:=r\transp Pr + 2\norm{B\transp P(F\mu_{w_t}+r)}_1 U_{\max}+ m\norm{B\transp PB}_\infty U_{\max}^2 +2|r\transp PF\mu_{w_t}|+\tr{F\transp PF\Sigma_{w_t}} $.
		 Since $x_t\transp A\transp P A x_t \le x_t\transp Px_t - x_t\transp x_t$, we have that
        \begin{align}
		\label{e:stab10}
        	\EE_{x_t}\bigl[x_{t+1}\transp Px_{t+1}\bigr] \leq x_t\transp Px_t - \norm{x_t}^2 + 2c_1\norm{x_t}_\infty + c_2.
        \end{align}
For $\theta\in\;]\max\{0,1-\lambda_{\max}(P)\}, 1[$ we know that
        \begin{align*}
         -\theta \norm{x_t}_\infty^2 + 2c_1\norm{x_t}_\infty + c_2 \le 0,\quad \forall \norm{x_t}_\infty > r,
        \end{align*}
		where $r := \frac{1}{\theta}\bigl(c_1 + \sqrt{c_1^2 + c_2\theta}\bigr)$. From~\eqref{e:stab10} it now follows that
		$
			\EE_{x_t}\bigl[x_{t+1}\transp P x_{t+1}\bigr] \le x_t\transp P x_t - (1-\theta)\norm{x_t}^2, \forall \norm{x_t}_\infty > r,
		$
		whence
		\[
			\EE_{x_t}\bigl[x_{t+1}\transp P x_{t+1}\bigr] \le \Bigl(1- \frac{1-\theta}{\lambdamax P}\Bigr)x_t\transp Px_t,\,\quad\forall\norm{x_t}_\infty > r.
		\]
		We see that the hypotheses of Lemma~\ref{l:FL} are satisfied with $V(x) := x\transp P x$, $\lambda := \left(1- \frac{1-\theta}{\lambdamax P}\right)$, and $K := \bigl\{x\in\R^n\big| \norm{x}_\infty \le r\bigr\}$. Since $\lambda_{\text{min}}(P)\norm{x}^2 \le x\transp Px$, it follows that
		\[
			\sup_{t\in\Nz} \EE_{\xz}\Bigl[\norm{x_t}^2\Bigr] \le \frac{1}{\lambda_{\text{min}}(P)}\sup\limits_{t\in\Nz}\EE_{\xz}\bigl[V(x_t)\bigr] < \infty,
		\]
which completes the proof.
	\end{proof}

\subsection{RHC Case}
In the RHC implementation is also iterative in nature, however instead of recalculating the gains at each time instant the optimization problem is solved every $kN$ steps, where $k\in\Nz$. The resulting optimal control policy (applied over a horizon $N$) is given by $\pi_{kN}^*=\bar G_{kN}^*\bar \varphi(\bar F\bar w)+\bar d_{kN}^*$, where again the control gains depend implicitly on the initial condition $x_{kN}$, i.e., $\bar G_{kN}^*=\bar G_{kN}^*(x_{kN}) $ and $\bar d_{kN}^*=\bar d_{kN}^*(x_{kN})$. Therefore, the optimal policy is given by
$\pi^{\rm RHC}=\left(\pi_{0}^*,\pi_{N}^*, \cdots\right)$. For $\ell = 1,\cdots, N$, the resulting closed-loop system over horizon $N$ is given by
\begin{equation}\label{eq:RHcl}
    x_{kN+\ell}=A^\ell x_{kN}+\bar B_\ell \bar G^*_{kN}\bar \varphi(\bar F\bar w)+\bar B_\ell \bar d^*_{kN}+\bar D_\ell\bar F\bar w + \bar D_\ell \bar r,
\end{equation}
where $k\in\Nz$, and  $\bar B_\ell$ and  $\bar D_\ell$ are suitably defined matrices that are extracted from $\bar B$ and $\bar D$, respectively.

\begin{proposition}\label{prop:RHmain}
	Assume that the matrix $A$ is Schur stable and the assumptions of Proposition \ref{p:main} hold. Then, under the control policy $\pi^{RHC}$ defined above, the closed loop system~\eqref{eq:RHcl} satisfies $\sup_{t\in\Nz} \EE_{\xz}\Bigl[\norm{x_t}^2\Bigr] < \infty$.
\end{proposition}
\begin{proof}
	Using \eqref{eq:RHcl} and the fact that $\EE_{x}\left[\bar\varphi(\bar F\bar w)\right]=0$, $\forall x\in\R^n$, we have that $\fa \ell = 1,\cdots, N$
	\begin{align*}
		& \EE_{x_{kN}}\bigl[x_{kN+\ell}\transp P_\ell x_{kN+l}\bigr]  =  x_{kN}\transp (A^\ell)\transp P_\ell A^\ell x_{kN}\\ & \, +2x_{kN}\transp (A^\ell)\transp P_\ell(\bar B_\ell \bar d^*_{kN}+\bar D_\ell\bar F\mu_{\bar w} + \bar D_\ell \bar r)+\bar r\transp\bar D_\ell\transp P_\ell \bar D_\ell \bar r\\&\, +(\bar d^*_{kN})\transp \bar B_\ell\transp P_\ell \bar B_\ell \bar d^*_{kN} +2(\bar d^*_{kN})\transp \bar B_\ell\transp P_\ell\bar D_\ell(\bar F\mu_{\bar w}+\bar r)\\&\, +2\mu_{\bar w}\transp \bar F\transp \bar D_\ell\transp P_\ell\bar D_\ell\bar r  + \tr{(\bar G^*_{kN})\transp \bar B_\ell \transp P_\ell \bar B_\ell \bar G_{kN}^*\Lambda_1}\\
		& \, +2\tr{(\bar G^*_{kN})\transp \bar B_\ell\transp P_\ell\bar D_\ell\bar F\Lambda_2} + \tr{\bar F\transp D_\ell\transp P_\ell D_\ell \bar F \Sigma_{\bar w}}.\nn
	\end{align*}
	Using the fact that $\norm{\bar d^*_{kN}}_\infty\leq U_{\rm max}$ and $\norm{\bar G_{kN}^*}_\infty \leq U_{\rm max}/\phi_{\rm max}$ (from ~\eqref{eq:problem2}), we obtain the following bound
	\begin{align*}
		& \EE_{x_{kN}}\bigl[x_{kN+\ell}\transp P_\ell x_{kN+\ell}\bigr] \\ &\quad \leq  x_{kN}\transp (A^\ell)\transp P_\ell A^\ell x_{kN} +2 c_{1\ell} \norm{x_{kN}}_\infty + c_{2\ell},
	\end{align*}
where $c_{1\ell}:=\norm{(A^\ell) \transp P_\ell \bar D_\ell (\bar F\mu_{\bar w}+\bar r)}_1+m\norm{(A^\ell)\transp P_\ell \bar B_\ell}_\infty U_{\rm max}$ and
$c_{2\ell}:=\bar r\transp \bar D_\ell\transp P\bar D_\ell \bar r + 2 \norm{(\bar B_\ell)\transp P_\ell \bar D_\ell (\bar F\mu_{\bar w}+\bar r)}_1 U_{\max}+ m\norm{\bar B_\ell\transp P_\ell \bar B_\ell}_\infty U_{\max}^2 +2|\bar r\transp \bar D_\ell\transp P_\ell \bar D_\ell \bar F\mu_{\bar w}|+\tr{\bar F\transp \bar D_\ell\transp P_\ell \bar D_\ell \bar F\Sigma_{\bar w}} +\max\limits_{\norm{\bar G^*_{kN}}_\infty\leq U_{\max}/\phi_{\max}}\big[\tr{\bar G^{*\mathsf T}_{kN} \bar B_\ell\transp P_\ell\bar B_\ell\bar G^*_{kN}\Lambda_1} +2\tr{\bar G^{*\mathsf T}_{kN}\bar B_\ell\transp P_\ell \bar D_\ell \bar F\Lambda_2}\big]$. Again, since $A$ is a Schur stable matrix (and hence $A^\ell$) there exists a matrix $P_\ell = P_\ell\transp>0$ with real valued entries that satisfies $(A^\ell)\transp P_\ell A^\ell - P_\ell \leq -\mathbf{I}_{n\times n}$, and its eigenvalues are real. Then we have $x_{kN}\transp(A^\ell)\transp P_\ell A^\ell x_{kN} \le x_{kN}\transp P_\ell x_{kN} - x_{kN}\transp x_{kN}$. Therefore,
	\begin{align}	
		\EE_{x_{kN}}\bigl[x_{kN+\ell}\transp P_\ell x_{kN+\ell}\bigr]& \leq x_{kN}\transp P_\ell x_{kN} - \norm{x_{kN}}^2\nn\\& + 2 c_{1\ell} \norm{x_{kN}}_\infty + c_{2\ell}.\label{e:stab11}
	\end{align}
	For $\theta_\ell\in\;]\max\{0,1-\lambda_{\max}(P_\ell)\}, 1[$ we know that
	\begin{align*}
		-\theta_\ell \norm{x_{kN}}_\infty^2 + 2 c_{1\ell}\norm{x_{kN}}_\infty +  c_{2\ell} \le 0,\,\forall\norm{x_{kN}}_\infty >  r_\ell,
	\end{align*}
	where $r_\ell := \frac{1}{\theta_\ell}\bigl( c_{1\ell}+ \sqrt{c_{1\ell}^2 + c_{2\ell}\theta_\ell}\bigr)$. From~\eqref{e:stab11} it now follows that
	$
		\EE_{x_{kN}}\bigl[x_{kN+\ell}\transp P_\ell x_{kN+\ell}\bigr] \le x_{kN}\transp P_\ell x_{kN} - (1-\theta_\ell)\norm{x_{kN}}^2, \forall \norm{x_{kN}}_\infty > r_\ell,
	$
	whence
	\begin{align}
	\label{eqn:boundl}
		\EE_{x_{kN}}\bigl[x_{kN+\ell}\transp P_\ell x_{kN+\ell}\bigr] \le \lambda_\ell x_{kN}\transp P_\ell x_{kN},\, \forall \norm{x_{kN}}_\infty > r_\ell,
	\end{align}
	where $\lambda_\ell := \Bigl(1- \frac{1-\theta}{\lambdamax {P_\ell}}\Bigr)$. Define $\lambda := \max\limits_{\ell=1,\cdots,N-1}\lambda_\ell$, $ r':=\max\limits_{\ell=1,\cdots,N-1}r_\ell$, $ \overline\lambda:= \max\limits_{\ell=1,\dots,N-1}\lambda_{\max}(P_\ell)$,  $\underline\lambda:= \min\limits_{\ell=1,\dots,N-1}\lambda_{\min}(P_\ell)$, then we can obtain using~\eqref{eqn:boundl} the conservative bound
	\begin{equation*}
	\label{eqn:boundl2}
		\EE_{x_{kN}}\bigl[x_{kN+\ell}\transp P_Nx_{kN+\ell}\bigr] \le \lambda' x_{kN}\transp P_Nx_{kN},\, \forall \norm{x_{kN}}_\infty > r'
	\end{equation*}
	for every $\ell = 1, \ldots, N-1$, where $\lambda':=\lambda\frac{\overline\lambda\lambda_{\max}(P_N)}{\underline\lambda \lambda_{\min}(P_N)}$, and the $N$-step bound
	\begin{align}
		\EE_{x_{kN}}\bigl[x_{(k+1)N}\transp P_Nx_{(k+1)N}\bigr] &\le \lambda_N x_{kN}\transp P_Nx_{kN},\nn \\
&\qquad  \forall \norm{x_{kN}}_\infty > r_N.\label{eqn:boundl23}
	\end{align}
	Let $V_N(x):=x\transp P_Nx$. Now, following the same reasoning as in Lemma~\ref{l:FL}, we can establish the following bound (for $k\in\Nz$, $\ell =1,\dots,N-1$)
	\begin{align}
		 &\EE_{x}\bigl[V_N(x_{kN+\ell})\bigr]=\EE_{x}\bigl[\EE[V_N(x_{kN+\ell})|x_{kN}]\bigr]\nn\\
		&\qquad\leq \EE_{x}\bigl[\EE[\lambda'V_N(x_{kN})+b'\indic{K'}(x_{kN})]\bigr]\nn\\
		&\qquad \leq \EE_{x}\bigl[\EE[\lambda'\EE[V_N(x_{kN})|x_{(k-1)N}]+b'\indic{K'}(x_{kN})]\bigr]\nn\\
		&\qquad \leq \EE_{x}\bigl[\EE[\lambda'\EE[\lambda_N V_N(x_{(k-1)N})+b\indic{K_N}(x_{(k-1)N})]\nn\\ &\qquad  \quad +b'\indic{K'}(x_{kN})]\bigr]\nn\\
		&\qquad \leq \lambda'\lambda_N^kV_N(x)+\sum_{i=0}^{k-1}\lambda_N^{k-1-i}b \EE_{x}\bigl[\indic{K_N}(x_{iN})\bigr]\nn\\ &\qquad \quad +b'\EE_x\bigl[\indic{K'}(x_{kN})\bigr]\nn\\
		&\qquad \leq \lambda'\lambda_N^kV_N(x)+\frac{b(1-\lambda_N^k)}{1-\lambda_N}+b', \label{eqn:finalbound}
	\end{align}
	where $b:= \sup\limits_{x\in K}\EE_x\bigl[V_N(x_{N})\bigr]$, $b':= \sup\limits_{x\in K'}\EE_x\bigl[V_N(x_{l})\bigr]$ for $\ell = 1,\cdots,N-1$, $K_N:=\bigl\{\xi\in\R^n\big|\norm{\xi}_\infty\leq r_N\bigr\}$, and $K':=\bigl\{\xi\in\R^n\big|\norm{\xi}_\infty\leq r'\bigr\}$. Note that the conditioning in the steps of~\eqref{eqn:finalbound} is done every $N$ steps as the problem is \emph{not Markovian} except then. Therefore, it follows from~\eqref{eqn:finalbound} that, $\fa t:=kN+\ell$,
	\begin{align}
		\sup\limits_{t\in\Nz}\EE_{x}\bigl[\norm{x_t}^2\bigr]&
		\leq \frac{1}{\lambda_{\min}(P_N)}\sup\limits_{t\in\Nz}\EE_{x}\bigl[V_N(x_{kN+l})\bigr],\nn\\
		&\hspace{-2cm}\leq \frac{1}{\lambda_{\min}(P_N)}\left(\lambda'\lambda_N^kV_N(x)+\frac{b}{1-\lambda_N}+b'\right)<\infty
	\end{align}
	which completes the proof.
\end{proof}

\subsection{Input-to-state Stability} \label{sec:iss}
	Input-to-state stability (\iss{}) is an interesting and important qualitative property of systems, dealing with input-output behavior. In the deterministic context~\cite{JiangWang-01} it generalizes the well-known bounded input bounded output (BIBO) property of linear systems~\cite[p.~490]{ref:antsaklisLS}. \iss{} provides a description of the behavior of a system subjected to bounded inputs. Here we are interested in a stochastic variant of input-to-state stability; see e.g.,~\cite{ref:borkarUnifStab, ref:tsiniasRobustStochISS} for other possible definitions and ideas (primarily in continuous-time).

	One possible way to measure the strength of stochastic inputs is in terms of their covariances; sometimes their moment generating functions are also employed. For Gaussian noise it is customary to consider a suitable norm of the covariance matrix as a measure of its strength. The deterministic version of input-to-state stability deals with $\mathcal L_\infty$-to-$\mathcal L_\infty$ gain from the input to the state of a system. We consider the linear system~\eqref{eq:system}, and establish a natural \iss{}-type property from the control and the noise inputs to the state of the system~\eqref{eq:system}, under both the MPC and the RHC strategies.

	\begin{defn}
		The system~\eqref{eq:system} is \emph{input-to-state stable in $\mathcal L_1$} if there exist functions $\beta\in\ClassKL$ and $\alpha, \gamma_1, \gamma_2\in\ClassKinfty$ such that for every initial condition $\xz\in\R^n$ and $\forall t\in\Nz$ we have
		\begin{equation}
		\label{eq:issm}
			\EE_{\xz}\bigl[\alpha(\norm{x_t})\bigr] \le \beta(\norm{\xz}, t) + \gamma_1\Bigl(\sup_{s\in\Nz}\norm{u_s}_\infty\Bigr) + \gamma_2\bigl(\norm{\Sigma}'\bigr),
		\end{equation}
		where $\norm{\cdot}'$ is an appropriate matrix norm.\DefEnd
	\end{defn}

	One difference with the deterministic definition of \iss{} is immediately evident, namely, the presence of the function $\alpha$ inside the expectation in~\eqref{eq:issm}. It turns out that often it is more natural to arrive at an estimate of $\EE_{\xz}\bigl[\alpha(\norm{x_t})\bigr]$ for some $\alpha\in\ClassKinfty$ than an estimate of $\EE_{\xz}[\norm{x_t}]$. Moreover, in case $\alpha$ is convex, Jensen's inequality~\cite[p.~348]{ref:dudley} shows that such an estimate implies an estimate of $\EE_{\xz}[\norm{x_t}]$. The following proposition can be easily established with the aid of Proposition~\ref{prop:SHmain} and Proposition~\ref{prop:RHmain}.

	\begin{proposition}
		The closed-loop systems~\eqref{eq:clsys} and~\eqref{eq:RHcl} are input-to-state stable in $\mathcal L_1$.\hfill $\blacksquare$
	\end{proposition}
The proof is omitted for space limitations.
\section{Numerical Example}\label{sec:Nexample}
\begin{figure}
\centering
	\subfigure[MPC implementation]
	{
  		\includegraphics[scale=0.35]{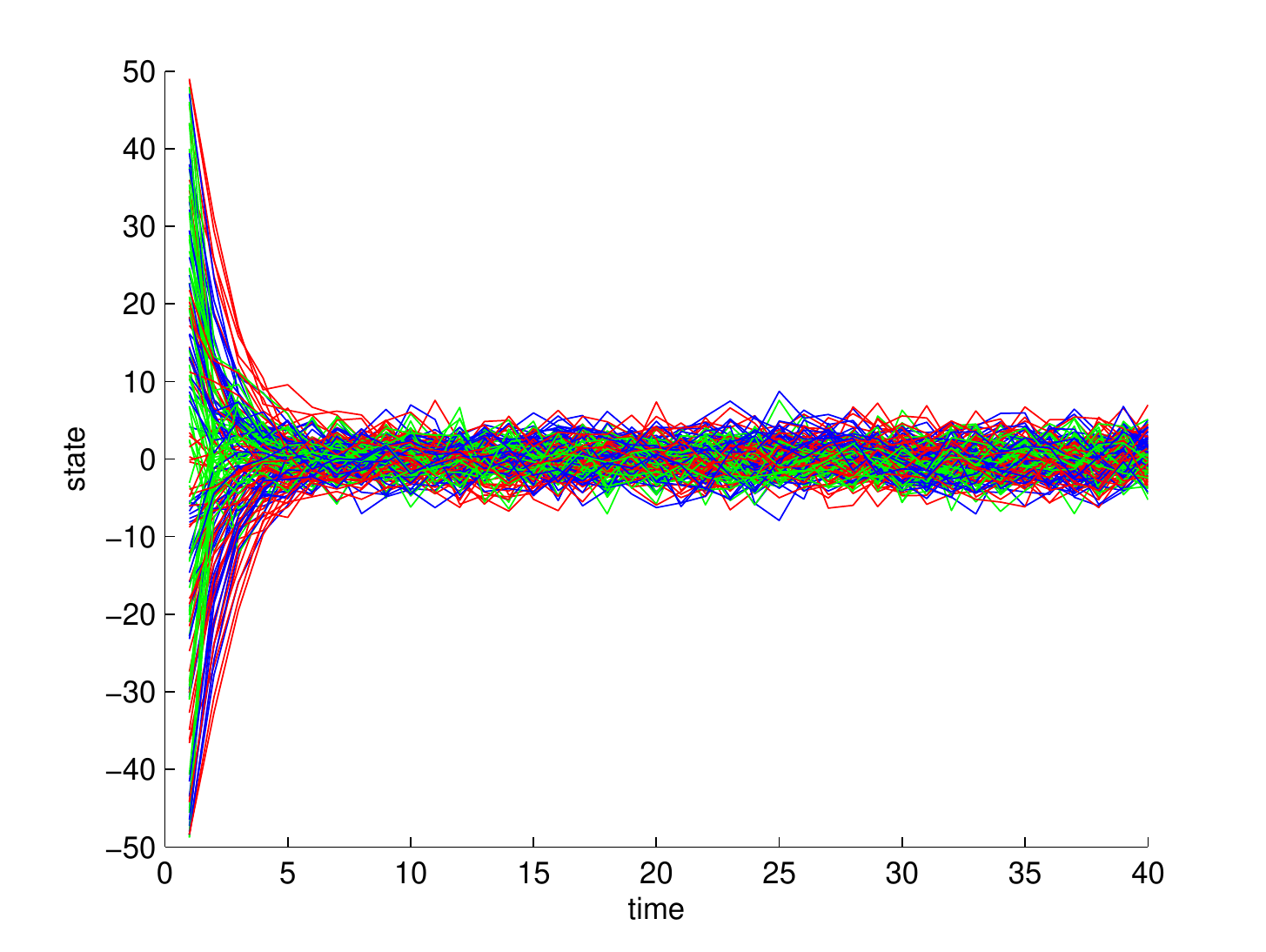}
		\label{fig:mpc}
	}
  	\subfigure[RHC implementation]
	{
		\includegraphics[scale=0.35]{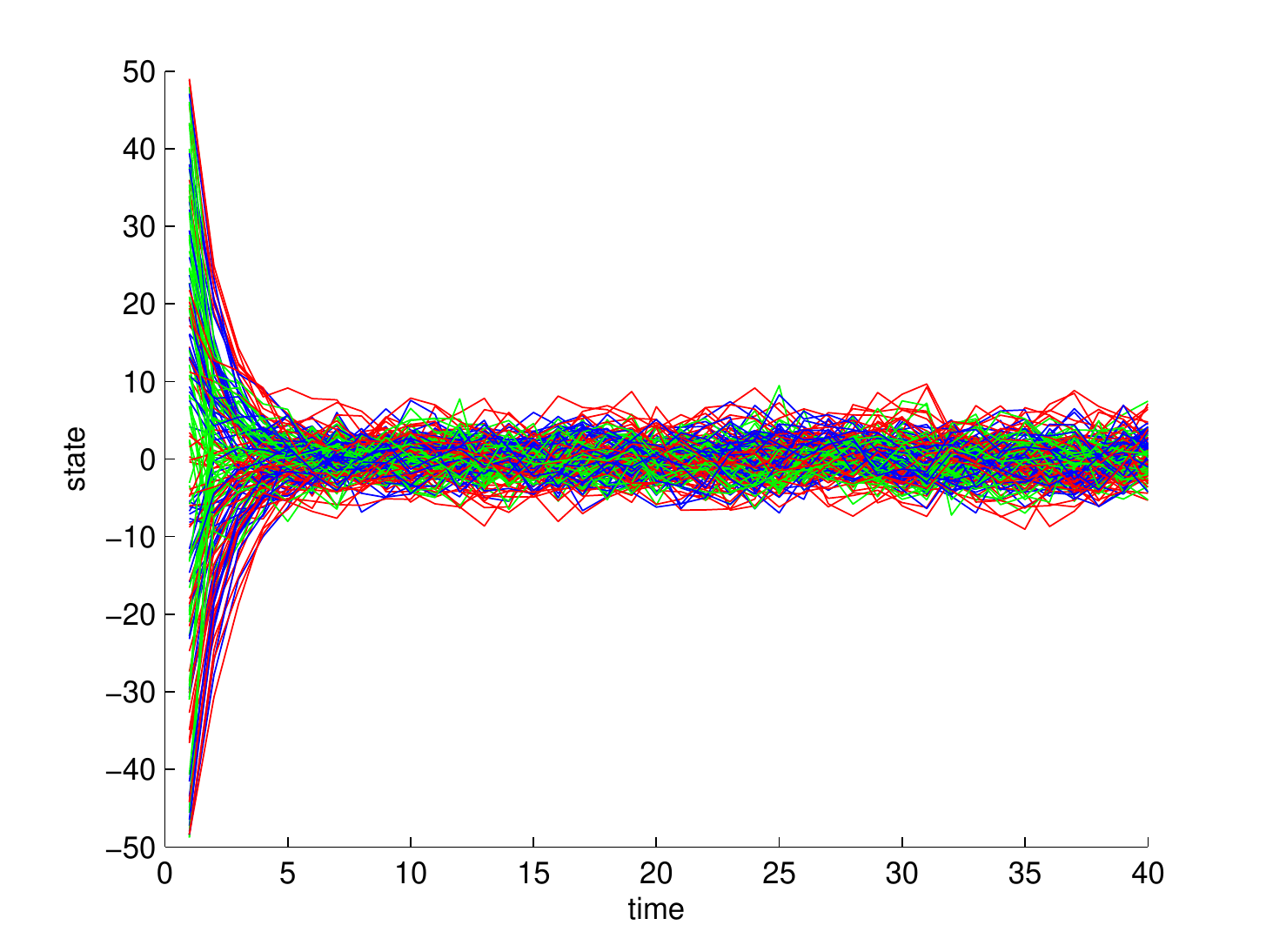}
		\label{fig:rh}
	}
	\caption{MPC and RHC algorithms corresponding to the system in~\S\ref{sec:Nexample}. The plots correspond to the aforementioned algorithms each run from 50 identical initial conditions distributed uniformly over  $[-50,50]$.}
\end{figure}

Let us consider the system \eqref{eq:system} with some generic  matrices $A=\smat{0.8 & 0.1 & 0.01\\ 0.3 & 0.3 & 0.06\\ 0.09 & 0.02 & 0.5}$, $B = \smat{1 \\ 2 \\ 0.5}$,
$F=\mathbf I_{3\times 3}$, and $r=\mathbf 0_{3\times 1}$.
We simulate the system starting from $50$ different initial conditions, all of which are sampled according to a uniform distribution over $[-50,50]^3$. The noise inputs are independent and identically sampled according to a normal distribution,  $w\sim \mathcal N(0,4\mathbf I_{3\times 3})$, the noise saturation function is chosen as in Example \ref{ex:sigmoids} with $\phi_{\max}=5$, and the input saturation bound $U_{\max}=10$. The optimization gain matrices are chosen to be $Q_i = 3\mathbf I_{3\times 3}$ and $R_i=2\mathbf I_{1\times 1}$, $\forall i$, and the optimization horizon $N=6$. The
optimization matrices are given by
$\Lambda_1 =  3.3024\mathbf I_{9\times 9}$, and $\Lambda_2 =  0.7846\mathbf I_{9\times 9}$. We used the \texttt{cvx} solver~\cite{ref:boydCVX} to handle the optimization problem \eqref{eq:problem2}. The results for the MPC implementation are shown in Figure \ref{fig:mpc}, and those for the RHC implementation are shown in Figure \ref{fig:rh}, for the full state evolution over a horizon of 40 time steps. Finally, it is interesting to note that the MPC and RHC \emph{average} performance indices over the 50 different runs are given by 3985 and 4327, respectively.

\section{Conclusions}\label{sec:conclusions}
In this paper, we provided a tractable optimization program that solves the stochastic Model Predictive Control and Rolling Horizon Control problems, while guaranteeing the satisfaction of hard bounds on the control input. We have showed that in both cases the resulting closed-loop process has bounded variance. We demonstrated that both implementations enjoy some qualitative notion of stochastic input-to-state stability. We provided several examples in which crucial matrices in our optimization program can be calculated off-line. Future direction for this research is aimed at lifting the current feedback strategy onto general vector spaces.


\begin{thebibliography}{10}
\providecommand{\url}[1]{#1}
\csname url@rmstyle\endcsname
\providecommand{\newblock}{\relax}
\providecommand{\bibinfo}[2]{#2}
\providecommand\BIBentrySTDinterwordspacing{\spaceskip=0pt\relax}
\providecommand\BIBentryALTinterwordstretchfactor{4}
\providecommand\BIBentryALTinterwordspacing{\spaceskip=\fontdimen2\font plus
\BIBentryALTinterwordstretchfactor\fontdimen3\font minus
  \fontdimen4\font\relax}
\providecommand\BIBforeignlanguage[2]{{%
\expandafter\ifx\csname l@#1\endcsname\relax
\typeout{** WARNING: IEEEtran.bst: No hyphenation pattern has been}%
\typeout{** loaded for the language `#1'. Using the pattern for}%
\typeout{** the default language instead.}%
\else
\language=\csname l@#1\endcsname
\fi
#2}}

\bibitem{ref:recstrat}
\BIBentryALTinterwordspacing
D.~Chatterjee, E.~Cinquemani, G.~Chaloulos, and J.~Lygeros, ``Stochastic
  optimal control up to a hitting time: optimality and rolling-horizon
  implementation,'' 2008. [Online]. Available:
  \url{http://arxiv.org/abs/0806.3008}
\BIBentrySTDinterwordspacing

\bibitem{MayneRawlingsRaoScokaert-00}
D.~Q. Mayne, J.~B. Rawlings, C.~V. Rao, and P.~O.~M. Scokaert, ``Constrained
  model predictive control: stability and optimality,'' \emph{Automatica},
  vol.~36, no.~6, pp. 789--814, Jun 2000.

\bibitem{BemporadMorari-99}
A.~Bemporad and M.~Morari, ``Robust model predictive control: A survey,''
  \emph{Robustness in Identification and Control}, vol. 245, pp. 207--226,
  1999.

\bibitem{ref:maciejowskibk}
J.~M. Maciejowski, \emph{Predictive {C}ontrol with {C}onstraints}.\hskip 1em
  plus 0.5em minus 0.4em\relax Prentice Hall, 2001.

\bibitem{ref:blanchini1999sic}
F.~Blanchini, ``{Set invariance in control},'' \emph{Automatica}, vol.~35,
  no.~11, pp. 1747--1767, 1999.

\bibitem{BertsimasBrown-07}
D.~Bertsimas and D.~B. Brown, ``Constrained stochastic {LQC}: a tractable
  approach,'' \emph{IEEE Transactions on Automatic Control}, vol.~52, no.~10,
  pp. 1826--1841, 2007.

\bibitem{Primbs-07}
J.~Primbs, ``A soft constraint approach to stochastic receding horizon
  control,'' in \emph{Proceedings of the 46th IEEE Conference on Decision and
  Control}, 2007, pp. 4797 -- 4802.

\bibitem{PrimbsSung-08}
J.~A. Primbs and C.~H. Sung, ``Stochastic receding horizon control of
  constrained linear systems with state and control multiplicative noise,''
  \emph{IEEE Trans. Automatic Control}, 2008, to appear.

\bibitem{ref:CannonKouvaritakisWu-08}
M.~Cannon, B.~Kouvaritakis, and X.~Wu, ``Probabilistic constrained {MPC} for
  systems with multiplicative and additive stochastic uncertainty,'' in
  \emph{IFAC World Congress}, Seoul, Korea, July 2008.

\bibitem{OldewurtelJonesMorari-08}
\BIBentryALTinterwordspacing
F.~Oldewurtel, C.~Jones, and M.~Morari, ``A tractable approximation of chance
  constrained stochastic {MPC} based on affine disturbance feedback,'' in
  \emph{Conference on Decision and Control, CDC}, Cancun, Mexico, Dec. 2008.
  [Online]. Available:
  \url{http://control.ee.ethz.ch/index.cgi?page=publications;action=details;id%
=3118}
\BIBentrySTDinterwordspacing

\bibitem{batinaPhDthesis}
I.~Batina, ``Model predictive control for stochastic systems by randomized
  algorithms,'' Ph.D. dissertation, Technische Universiteit Eindhoven, 2004.

\bibitem{MaciejowskiLecchiniLygeros-05}
M.~Maciejowski, A.~Lecchini, and J.~Lygeros, ``{NMPC} for complex stochastic
  systems using {Markov Chain Monte Carlo},'' in \emph{International Workshop
  on Assessment and Future Directions of Nonlinear Model Predictive Control},
  ser. Lecture Notes in Control and Information Sciences, vol. 358/2007.\hskip
  1em plus 0.5em minus 0.4em\relax Stuttgart, Germany: Springer, 2005, pp.
  269--281.

\bibitem{ref:Stoorvogel}
\BIBentryALTinterwordspacing
A.~A. Stoorvogel, A.~Saberi, and S.~Weiland, ``On external semi-global
  stochastic stabilization of linear systems with input saturation,'' 2006,
  submitted. [Online]. Available:
  \url{http://homepage.mac.com/a.a.stoorvogel/subm03.pdf}
\BIBentrySTDinterwordspacing

\bibitem{AgarwalCinquemaniChatterjeeLygeros-09}
\BIBentryALTinterwordspacing
M.~Agarwal, E.~Cinquemani, D.~Chatterjee, and J.~Lygeros, ``On convexity of
  stochastic optimization problems with constraints,'' in \emph{European
  Control Conference}, 2009, submitted. [Online]. Available:
  \url{http://control.ee.ethz.ch/index.cgi?page=publications;action=details;id%
=3271}
\BIBentrySTDinterwordspacing

\bibitem{AldenSmith-92}
J.~M. Alden and R.~L. Smith, ``Rolling horizon procedures in nonhomogeneous
  {M}arkov decision processes,'' \emph{Operations Research}, vol.~40, no.
  suppl. 2, pp. S183--S194, May-Jun. 1992.

\bibitem{ref:YangSontagSussmann-97}
Y.~D. Yang, E.~D. Sontag, and H.~J. Sussmann, ``Global stabilization of linear
  discrete-time systems with bounded feedback,'' \emph{Systems and Control
  Letters}, vol.~30, no.~5, pp. 273--281, 1997.

\bibitem{ref:ben-tal04}
A.~Ben-Tal, A.~Goryashko, E.~Guslitzer, and A.~Nemirovski, ``Adjustable robust
  solutions of uncertain linear programs,'' \emph{Mathematical Programming},
  vol.~99, no.~2, pp. 351--376, 2004.

\bibitem{ref:goulart06}
P.~J. Goulart, E.~C. Kerrigan, and J.~M. Maciejowski, ``Optimization over state
  feedback policies for robust control with constraints,'' \emph{Automatica J.
  IFAC}, vol.~42, no.~4, pp. 523--533, 2006.

\bibitem{Luenberger-69}
D.~Luenberger, \emph{Optimization by Vector Space Methods}.\hskip 1em plus
  0.5em minus 0.4em\relax J. Wiley \& Sons, 1969.

\bibitem{ref:boyd04}
S.~Boyd and L.~Vandenberghe, \emph{Convex {O}ptimization}.\hskip 1em plus 0.5em
  minus 0.4em\relax Cambridge: Cambridge University Press, 2004, sixth printing
  with corrections, 2008.

\bibitem{ref:boydCVX}
M.~Grant and S.~Boyd, ``{CVX}: {M}atlab software for disciplined convex
  programming (web page and software),'' \url{http://stanford.edu/~boyd/cvx},
  December 2000.

\bibitem{ref:AbramowitzStegun}
M.~Abramowitz and I.~A. Stegun, \emph{Handbook of Mathematical Functions with
  Formulas, Graphs, and Mathematical Tables}, ser. National Bureau of Standards
  Applied Mathematics Series.\hskip 1em plus 0.5em minus 0.4em\relax For sale
  by the Superintendent of Documents, U.S. Government Printing Office,
  Washington, D.C., 1964, vol.~55.

\bibitem{ref:meyn93}
S.~P. Meyn and R.~L. Tweedie, \emph{Markov {C}hains and {S}tochastic
  {S}tability}.\hskip 1em plus 0.5em minus 0.4em\relax London: Springer-Verlag,
  1993.

\bibitem{ref:foss04}
S.~Foss and T.~Konstantopoulos, ``An overview of some stochastic stability
  methods,'' \emph{Journal of Operations Research Society of Japan}, vol.~47,
  no.~4, pp. 275--303, 2004.

\bibitem{JiangWang-01}
Z.-P. Jiang and Y.~Wang, ``Input-to-state stability for discrete-time nonlinear
  systems,'' \emph{Automatica}, vol.~37, no.~6, pp. 857--869, June 2001.

\bibitem{ref:antsaklisLS}
P.~J. Antsaklis and A.~N. Michel, \emph{Linear {S}ystems}.\hskip 1em plus 0.5em
  minus 0.4em\relax Boston, MA: Birkh\"auser Boston Inc., 2006.

\bibitem{ref:borkarUnifStab}
V.~S. Borkar, ``Uniform stability of controlled {M}arkov processes,'' in
  \emph{System theory: modeling, analysis and control (Cambridge, MA, 1999)},
  ser. Kluwer International Series in Engineering Computer Science.\hskip 1em
  plus 0.5em minus 0.4em\relax Boston, MA: Kluwer Academic Publishers, 2000,
  vol. 518, pp. 107--120.

\bibitem{ref:tsiniasRobustStochISS}
J.~Spiliotis and J.~Tsinias, ``Notions of exponential robust stochastic
  stability, {ISS} and their {L}yapunov characterization,'' \emph{International
  Journal of Robust and Nonlinear Control}, vol.~13, no.~2, pp. 173--187, 2003.

\bibitem{ref:dudley}
R.~M. Dudley, \emph{Real {A}nalysis and {P}robability}, ser. Cambridge Studies
  in Advanced Mathematics.\hskip 1em plus 0.5em minus 0.4em\relax Cambridge:
  Cambridge University Press, 2002, vol.~74, revised reprint of the 1989
  original.


\end{thebibliography}
\end{document}